\documentclass[psamsfonts,reqno]{amsart}
\usepackage{amssymb,eucal,latexsym,xy}
\xyoption{matrix}\xyoption{arrow}\xyoption{ps}
\theoremstyle{plain}

\newtheorem{lemma}{Lemma}[section]
\newtheorem{theorem}[lemma]{Theorem}
\newtheorem{proposition}[lemma]{Proposition}
\newtheorem{corollary}[lemma]{Corollary}

\newtheorem*{stat}{\name}
\newcommand{\name}{testing}

\theoremstyle{definition}
\newtheorem*{definition}{Definition}

\theoremstyle{remark}
\newtheorem{remark}[lemma]{Remark}

\newcommand{\qedc}{{\qed}~{\rm Claim~{\theclaim}.}}
\newcommand{\qedsc}{{\qed}~{\rm Claim.}}

\numberwithin{equation}{section}

\newcommand{\two}{\mathbf{2}}
\newcommand{\eps}{\varepsilon}

\newcommand{\id}{\mathrm{id}}

\newcommand{\tvi}{\vrule height 12pt depth 6pt width 0pt}

\newcommand{\di}[1]{\epsilon({#1})}
\newcommand{\DI}[1]{\Lambda(#1)}

\newcommand{\AC}[2]{G_{#1}[#2]}

\newcommand{\fm}{\mathfrak{m}}
\newcommand{\fn}{\mathfrak{n}}

\newcommand{\pup}[1]{\textup{(}{#1}\textup{)}}

\newcommand{\NN}{\mathbb{N}}
\newcommand{\ZZ}{\mathbb{Z}}

\newcommand{\zo}{0}
\newcommand{\nz}[1]{{#1}^{\sqcup0}}
\newcommand{\into}{\hookrightarrow}

\newcommand{\znzo}{(\ZZ/n\ZZ)\sqcup\set{0}}
\newcommand{\znznz}{(\ZZ/n\ZZ)^{\sqcup0}}

\newcommand{\jirr}{join-ir\-re\-duc\-i\-ble}

\newcommand{\jirry}{join-ir\-re\-duc\-i\-bil\-i\-ty}

\newcommand{\ram}{regular abel\-ian mon\-oid}
\newcommand{\jz}{$(\vee,0)$}

\newcommand{\jzs}{\jz-semi\-lat\-tice}

\newcommand{\set}[1]{\left\{#1\right\}}
\newcommand{\setm}[2]{\set{{#1}\mid{#2}}}
\newcommand{\famm}[2]{\left({#1}\mid{#2}\right)}

\DeclareMathOperator{\J}{J} \DeclareMathOperator{\Id}{Id}
\DeclareMathOperator{\Sub}{Sub}

\newcommand{\refmatr}[8]{
  \begin{tabular}{|c|c|c|}
  \cline{2-3}
  \multicolumn{1}{l|}{} & ${#1}$ & ${#2}$\tvi\\
  \hline
  ${#3}$ & ${#4}$ & ${#5}$\tvi\\
  \hline
  ${#6}$ & ${#7}$ & ${#8}$\tvi\\
  \hline
  \end{tabular}
  }

\newcommand{\bD}{\boldsymbol{D}}

\newcommand{\cB}{\mathcal{B}}

\newcommand{\cF}{\mathcal{F}}
\newcommand{\cL}{\mathcal{L}}
\newcommand{\Rep}{\mathcal{R}_{\mathrm{ep}}}
\newcommand{\cO}{\mathcal{O}}

\newcommand{\mbar}{\overline{m}}

\begin{document}

\author{K.\,R. Goodearl}
\address{Department of Mathematics\\
University of California\\
Santa Barbara, CA 93106\\
U.S.A.} \email{goodearl@math.ucsb.edu}
\urladdr{http://www.math.ucsb.edu/\~{}goodearl/}

\author{E. Pardo}
\address{Departamento de Matem\'aticas\\
Universidad de C\'adiz\\
Apartado 40\\
11510 Puerto Real (C\'adiz)\\
Spain} \email{enrique.pardo@uca.es}
\urladdr{http://matematicas.uca.es/\~{}neumann/EMAIN.HTML}

\author{F. Wehrung}
\address{CNRS, UMR 6139\\
D\'epartement de Math\'ematiques\\
Universit\'e de Caen\\
14032 Caen Cedex\\
France} \email{wehrung@math.unicaen.fr}
\urladdr{http://www.math.unicaen.fr/\~{}wehrung}

\title[Semilattices of groups and limits of Cuntz
algebras]{Semilattices of groups and inductive limits\\
of Cuntz algebras}

\thanks{The research of the first author was partially supported by
an NSF grant. The research of the second author was partially
supported by DGESIC grant BFM2001--3141, by FQM-298 PAI III grant
(Junta de Andaluc\'\i a), and by PR2001-0276 grant from Secretar\'\i a de
Estado de Educaci\'on y Universidades, M.E.C.D. (Spain). The
research of the third author was partially supported by INTAS
project 03-51-4110}

\subjclass[2000]{Primary 20M17, 46L35; Secondary 06A12, 06F05,
46L80} \keywords{Cuntz algebra, Cuntz limit, direct limit,
inductive limit, non-stable K-theory, regular monoid, Riesz
refinement, Mayer-Vietoris property, strongly periodic,
semilattice, abelian group, pure subgroup}

\begin{abstract}
We characterize, in terms of elementary properties, the abelian
monoids which are direct limits of finite direct sums of monoids
of the form $\znzo$ (where $\zo$ is a new zero element), for
positive integers $n$. The key properties are the Riesz refinement
property and the requirement that each element $x$ has finite
order, that is, $(n+1)x=x$ for some positive integer~$n$. Such
monoids are necessarily semilattices of abelian groups, and part
of our approach yields a characterization of the Riesz refinement
property among semilattices of abelian groups. Further, we describe
the monoids in question as certain submonoids of direct products
$\Lambda\times G$ for semilattices $\Lambda$ and torsion abelian
groups $G$. When applied to the monoids
$V(A)$ appearing in the non-stable K-theory of C*-algebras, our
results yield characterizations of the monoids
$V(A)$ for C* inductive limits $A$ of sequences of finite direct
products of matrix algebras over Cuntz algebras $\cO_n$. In
particular, this completely solves the problem of determining the
range of the invariant in the unital case of R\o rdam's
classification of inductive limits of the above type.
\end{abstract}

\maketitle
\section{Introduction}\label{S:Intr}

As indicated in the abstract, the goal of this paper is to prove a
semigroup-the\-or\-et\-ic result motivated by, and with
applications to, the classification theory of C*-algebras. The
relevant C*-algebras, which we will call \emph{Cuntz limits\/} for
short, are the C* inductive limits of sequences of finite direct
products of full matrix algebras over the Cuntz algebras $\cO_n$.
(We recall the definition of the latter for the information of
non-C*-algebraic readers: for $2\leq n<\infty$, the Cuntz algebra
$\cO_n$, introduced in~\cite{Cuntz}, is the unital C*-algebra
generated by elements $s_1$,\dots, $s_n$ with relations
$s_i^*s_j=\delta_{ij}$ and $\sum_{i=1}^n s_is_i^*=1$.) Our results
will provide an analogue for Cuntz limits of the description of
the range of the invariant for separable AF C*-algebras (namely,
ordered $K_0$) by Elliott~\cite{Elli} and Effros, Handelman, and
Shen~\cite{EHS}. We begin by sketching the source of the problem
and giving a precise formulation. Most of the remainder of the
paper is purely semigroup-theoretic, except for the applications
to C*-algebras in the final section.

In \cite{Ro}, R\o rdam gave a K-theoretic classification of even
Cuntz limits (i.e., C* inductive limits of sequences of finite
direct products of matrix algebras over $\cO_n$s with $n$ even).
The invariant which R\o rdam used for his classification is
equivalent, in the unital case, to the pair $(V(A),[1_A])$ where
$V(A)$ denotes the (additive, commutative) monoid of Murray-von
Neumann equivalence classes of projections (self-adjoint
idempotents) in matrix algebras over a C*-al\-ge\-bra~$A$, and
$[1_A]$ is the class in $V(A)$ of the unit projection in $A$
(cf.~\cite[Sections 4.6, 5.1, and 5.2]{Black}). Thus, the unital
case of the classification states that if $A$ and $B$ are unital
even Cuntz limits, then $A\cong B$ if and only if $(V(A),[1_A])
\cong (V(B),[1_B])$, that is, there is a monoid isomorphism $V(A)
\rightarrow V(B)$ sending $[1_A]$ to $[1_B]$
(cf.~\cite[Theorem~7.1]{Ro}). R\o rdam has communicated to us
\cite{Rotwo} that his classification can be extended to all Cuntz
limits by investing the work of Kirchberg \cite{Kirch} and
Phillips~\cite{Phil}.

As with any classification theorem, an accompanying problem is to
describe the range of the invariant -- that is, which pairs
$(M,u)$ (an abelian monoid $M$ together with an element $u\in M$)
appear as $(V(A),[1_A])$ for unital  Cuntz limits $A$? This
question reduces to an interesting problem in the theory of
monoids which we shall describe shortly. The major aim of this
paper is to solve this monoid problem, and then draw corresponding
conclusions for Cuntz limits. For non-unital Cuntz limits~$A$, R\o
rdam's classifying invariant amounts to a triple
$(V(A),P(A),\tau)$ where $P(A)$ is a partial semigroup consisting
of unitary equivalence classes of projections in $A$ and
$\tau\colon P(A)\rightarrow V(A)$ is a natural homomorphism. Thus,
$V(A)$ is an important part of the classification in general, and
pinning down its structure is of interest also in the non-unital
case.

In trying to match a given pair $(M,u)$ with a unital Cuntz limit,
it is easy to eliminate~$u$. First, one notes that $u$ must be an
\emph{order-unit\/} in $M$, that is, for any $x\in M$, there exist
$y\in M$ and $n\in\NN$ such that $x+y= nu$. Second, if we can find
a Cuntz limit~$B$ such that $V(B)\cong M$, then there is a
projection $p$ in some matrix algebra $M_n(B)$ whose class $[p]$
corresponds to $u$, and the C*-algebra $A=pM_n(B)p$ is a unital
Cuntz limit satisfying $(V(A),[1_A])\cong (M,u)$. Thus, we
concentrate on the problem of describing those abelian monoids
which appear as $V(A)$s. In the case of simple algebras, R\o
rdam's work provides the answer -- the abelian monoids appearing
as $V(A)$ for simple (unital) Cuntz limits $A$ are the monoids
$G\sqcup \set{\zo}$ for arbitrary countable torsion abelian groups
$G$ \cite[Proposition~2.5 and Theorem~2.6]{Ro}, where $G\sqcup
\set{\zo}$ is the monoid obtained from~$G$ by adjoining a new zero
element. The answer is also known for the case of $\cO_2$-limits
(Cuntz limits involving only direct products of matrix algebras
over $\cO_2$), one of the basic ingredients of a class of
C*-algebras classified by Lin in \cite{Li}. The monoids appearing
as $V(A)$ for $\cO_2$-limits are just the direct limits of
sequences of \emph{Boolean\/} monoids (finite direct sums of
copies of the two-element monoid). These direct limits were shown
by Bulman-Fleming and McDowell to be precisely the countable
distributive upper semilattices, see \cite[Theorem~3.1]{BFMD}.
While the result of \cite{BFMD} relies heavily on Shannon's
categorical result \cite[Theorem~2]{Shan}, a purely general
algebraic proof has been given by the first and third authors
\cite[Theorem~6.6]{GW}.

It is known that the functor $V(-)$ converts C* inductive limits
to monoid inductive (direct) limits, that it converts finite
direct products to direct sums, and that $V(M_m(A))\cong V(A)$ for
all $A$ and $m$. Moreover,
$V(\cO_n)\cong(\ZZ/(n-1)\ZZ)\sqcup\set{\zo}$ (this follows from
the computations in \cite{Cuntz2}; see also
Section~\ref{S:C*applics}). Thus, the monoid problem boils down to
the following task (where we have replaced $n-1$ by $n$ for
convenience):
  \begin{quote}\em
  Characterize those abelian monoids isomorphic to
  direct limits of sequences of finite direct sums of building
  blocks of the form $\znzo$.
  \end{quote}
In this paper, we solve the above problem, and thus characterize
the monoids that appear as $V(A)$ for Cuntz limits $A$.

\section{Background}\label{S:Back}

\subsection*{Monoids} All monoids in this paper will be abelian,
written additively, and so with additive identities denoted $0$.
The monoids that appear as $V(A)$ for Cuntz limits $A$ enjoy
several standard properties familiar from other classification
results, such as conicality and refinement. Recall that a monoid
$M$ is \emph{conical\/} if $x+y=0$ (for $x$, $y\in M$) always
implies $x=y=0$, and that $M$ satisfies the \emph{Riesz refinement
property\/} provided that for any $x_1,x_2,y_1,y_2\in M$
satisfying $x_1+x_2= y_1+y_2$, there exist elements $z_{ij}\in M$
such that each $x_i= z_{i1}+z_{i2}$ and each $y_j= z_{1j}+z_{2j}$.
It is convenient to record the latter four equations in the form
of a \emph{refinement matrix\/}:
  \[
  \refmatr{y_1}{y_2}{x_1}{z_{11}}{z_{12}}{x_2}{z_{21}}{z_{22}}
  \]
Following \cite{Dobb}, a \emph{refinement monoid\/} is any abelian
monoid satisfying the Riesz refinement property.

Any abelian monoid $M$ supports a translation-invariant pre-order
$\leq$ (often called the \emph{algebraic pre-order\/}) defined by
the existence of differences: $x\leq y$ if and only if there
exists $z\in M$ such that $x+z=y$. All inequalities in abelian
monoids will be with respect to this pre-order. The monoid $M$
satisfies the \emph{Riesz decomposition property\/} provided that
whenever $x\leq y_1+y_2$ in $M$, there exist elements $x_i\in M$
such that $x= x_1+x_2$ and each $x_i\leq y_i$. This property
follows from the refinement property, but in general the two are
not equivalent.

We can construct a monoid from any additive group $G$ by adjoining
a new additive identity, denoted $0$ following our general
convention. The new monoid can be expressed in the form
$G\sqcup\set{0}$, which we sometimes abbreviate $\nz{G}$. In case
we need to refer to the zero of the group $G$, we write $0_G$ in
order to distinguish this element from the zero of the monoid
$\nz{G}$.

Let $M$ be an abelian monoid and $x\in M$. It is standard in the
semigroup literature to say that $x$ is \emph{periodic\/} if the
subsemigroup of $M$ generated by $x$ is finite. This does not,
however, imply that this subsemigroup is a group. Thus, we shall
say that $x$ is \emph{strongly periodic\/} provided the
subsemigroup generated by~$x$ is a finite group; note that this
occurs if and only if there is a positive integer~$m$ such that
$(m+1)x=x$. The smallest such $m$ is, of course, the order of the
sub(semi)group generated by $x$; we will refer to it as the
\emph{order\/} of $x$. We say that~$M$ itself is \emph{strongly
periodic\/} provided every element of $M$ is strongly periodic.
\medskip

\subsection*{Semilattices} Recall that an \emph{upper
semilattice\/} (or \emph{$\vee$-semilattice\/}) is a partially
ordered set in which every pair of elements has a supremum. All
semilattices in this paper will be upper semilattices, and they
will also be assumed to have least elements, denoted $0$. We will
refer to them simply as \emph{semilattices\/}, rather than using
the precise but cumbersome term ``\jzs''. If one takes $+=\vee$,
any semilattice becomes an abelian monoid in which $2x=x$ for all
$x$; conversely, any abelian monoid with the latter property is a
semilattice with respect to its algebraic pre-order. (It is an
easy exercise to check that the pre-order is actually a partial
order in this case.) Thus, for our purposes, it is convenient to
take the name ``semilattice'' to mean any abelian monoid in which
all elements satisfy the equation $2x=x$. Note that in a
semilattice, $x\leq y$ if and only if $x+y=y$. We shall generally
write the operation in a semilattice as addition, except when it
appears helpful to emphasize that an element $x\vee y$ is the
supremum of elements $x$ and $y$.

An \emph{ideal\/} of a semilattice $S$ is any nonempty,
order-hereditary subset which is closed under finite suprema, that
is, any submonoid of $S$ which is hereditary with respect to the
algebraic order. The collection of ideals of $S$ is a complete
lattice, denoted $\Id S$, in which infima are given by
intersections. There is a canonical semilattice embedding of $S$
into $\Id S$ given by $a\mapsto [0,a]$, where $[0,a]$ denotes the
``closed interval'' $\setm{x\in S}{x\leq a}$.

A \emph{distributive semilattice\/} is any semilattice which
satisfies the Riesz decomposition (equivalently, refinement)
property (cf. \cite[Lemma~2.3]{GW}). A semilattice $S$ is
distributive if and only if the ideal lattice $\Id S$ is
distributive \cite[Section~II.5]{GLT2}.
\medskip

\subsection*{Semilattices of Groups} Let $M$ be an abelian
monoid, and let $\DI{M}$ denote the set of idempotent (actually
``idem-multiple'') elements of $M$, that is, those $e\in M$ such
that $2e=e$. Then $\DI{M}$ is a submonoid of $M$, and it is a
semilattice. Note that the algebraic (pre-) order within $\DI{M}$
coincides with the restriction of the pre-order from $M$: if $e$,
$f\in \DI{M}$ and $e\leq f$ in $M$, then $e+x=f$ for some $x\in
M$, whence $e+f= 2e+x= e+x= f$, and so $e\leq f$ within $\DI{M}$.
Consequently, we may use inequalities for idempotents with no
danger of ambiguity.

The monoid $M$ is a \emph{semilattice of groups\/} provided $M$ is
a disjoint union of subgroups, that is, a disjoint union of
subsemigroups each of which happens to be a group. (The collection
of these subgroups is then a semilattice, where the supremum of
subgroups $G$ and $G'$ is the unique subgroup containing $G+G'$.)
The zero elements of these groups are then the idempotent elements
of $M$, and so $M$ will be a disjoint union of subgroups
$\AC{M}{e}$ indexed by the idempotents $e\in \DI{M}$. These
subgroups may be described as follows:
  \[
  \AC{M}{e}= \setm{x\in M}{e\leq x\leq e}.
  \]
Note that whenever $e\leq f$ in $\DI{M}$, the rule $x\mapsto x+f$
defines a group homomorphism $\AC{M}{e} \rightarrow \AC{M}{f}$.

If $M$ is a semilattice of groups, then the homomorphisms above,
together with the groups $\AC{M}{e}$, define a functor from
$\DI{M}$ (made into a category from its poset structure in the
standard way) to the category of abelian groups. Conversely (e.g.,
\cite[Theorem 4.11]{CP} or \cite[p.~89--90]{Ho}), given any
functor $\cF$ from a semilattice~$\Lambda$ to abelian groups, we
can construct a corresponding semilattice of groups, say
$M(\Lambda,\cF)$, whose underlying set is the disjoint union of
the groups $\cF(e)$ for $e\in \Lambda$. The addition operation in
$M(\Lambda,\cF)$ is defined as follows: if $x$, $y\in
M(\Lambda,\cF)$, there are unique $e$, $f\in \Lambda$ such that
$x\in \cF(e)$ and $y\in \cF(f)$, and $x+y := \cF(i)(x)+\cF(j)(y)$
in $\cF(e+f)$, where $i\colon e\rightarrow e+f$ and $j\colon
f\rightarrow e+f$ are the unique morphisms in the category
$\Lambda$ corresponding to the relations $e\leq e+f$ and $f\leq
e+f$.

Semilattices of groups are characterized by the standard
semigroup-theoretic concept of regularity, which takes the
following form in additive notation. An abelian monoid $M$ is
\emph{\pup{von Neumann} regular\/} provided that for each $x\in
M$, there exists $y\in M$ such that $x+y+x=x$. Equivalently, $M$
is regular if and only if $2x\leq x$ for all $x\in M$. Observe
that every strongly periodic monoid is regular.

It is well known that a semigroup $S$ (not necessarily
commutative) is a semilattice of groups if and only if $S$ is
regular and its idempotents are central \cite[Theorem~2.1]{Ho}. We
give a short proof of the commutative case below, for the reader's
convenience.

\begin{lemma}\label{L:SemGrps}
An abelian monoid $M$ is a semilattice of groups if and only if
$M$ is regular.
\end{lemma}

\begin{proof} $(\Longrightarrow)$: Any $x\in M$ lies in a group
$\AC{M}{e}$, for some $e\in \DI{M}$. Then $x+y=e$ for some $y\in
\AC{M}{e}$, whence $2x+y= x$.

$(\Longleftarrow)$: For $e\in \DI{M}$, set $X(e)= \setm{x\in
M}{e\leq x\leq e}$, and observe that $X(e)$ is a subsemigroup of
$M$, containing $e$. If $x\in X(e)$, there exist $y,z\in M$ such
that $e+y=x$ and $x+z=e$. Then $e+x= 2e+y= e+y= x$, which shows
that $e$ is an additive identity for $X(e)$. Since $z\leq e$, we
see that $z+e \in X(e)$, and then since $x+(z+e)=2e=e$, we see
that $z+e$ is an additive inverse for $x$ within $X(e)$. Therefore
$X(e)$ is a group.

It remains to prove that $M$ is the disjoint union of the groups
$X(e)$. Disjointness is clear, since if $x\in X(e)\cap X(f)$ for
some $e$, $f\in \DI{M}$, then $e\leq x\leq f\leq x\leq e$, whence
$e=f$. Given $x\in M$, we have $2x\leq x$ by hypothesis, whence
$2x+y=x$ for some $y\in M$. Set $e= x+y$, observing that $e\leq
x\leq e$ and $2e= 2x+y+y= x+y=e$, that is, $e\in \DI{M}$ and $x\in
X(e)$. Therefore $M$ is the disjoint union of the subgroups
$X(e)$, as desired.
\end{proof}

In view of Lemma~\ref{L:SemGrps}, the terms ``semilattice of
abelian groups'' and ``\ram'' are equivalent; we shall use the
latter from now on.

If $M$ is a \ram, then each element $a\in M$ lies in a group
$\AC{M}{\di{a}}$ for a unique idempotent $\di{a}\in \DI{M}$. Let
$a^-$ denote the additive inverse of~$a$ in the group
$\AC{M}{\di{a}}$.

\section{Regular refinement monoids}
\label{S:PurInfRefm}

We begin by establishing some necessary conditions for the general
type of direct limits that we are seeking to characterize, among
which are the key properties of regularity and refinement. We also
develop a new characterization of regular refinement monoids.

\begin{proposition}\label{P:NecessLL}
Let $M$ be any direct limit of finite direct sums of monoids of
the form $\nz{A}$, for abelian groups $A$. Then the following
statements hold:
\begin{itemize}
\item[(a)] $M$ is a regular conical refinement monoid.

\item[(b)] If all the groups $A$ are torsion groups, then $M$ is
strongly periodic.

\item[(c)] For any idempotents $e\leq f$ in $M$, the homomorphism
$x\mapsto x+f$ from $\AC{M}{e}$ to $\AC{M}{f}$ is injective.

\item[(d)] For any idempotents $e\leq f$ in $M$, the group
$\AC{M}{e}+f$ is a pure subgroup of $\AC{M}{f}$.
\end{itemize}
\end{proposition}

\begin{proof} Statement (b) is clear.
Note that (c) and (d) are equivalent to the following properties:
  \begin{itemize}
  \item[(c$'$)] If $e\leq f$ in $\DI{M}$ and $x\in \AC{M}{e}$ such
that $x+f=f$, then $x=e$.
  \item[(d$'$)] If $e\leq f$ in
$\DI{M}$ and $x\in \AC{M}{e}$, $y\in \AC{M}{f}$ satisfy $x+f=my$
for some $m\in\NN$, then there exists $z\in \AC{M}{e}$ such that
$x+f= m(z+f)$.
  \end{itemize}
Thus, properties (a), (c), (d) can all be checked in terms of
finite sets of equations involving finitely many elements.
Therefore we need only verify them in the case when $M=\nz{A}$.

(a) Obviously $M$ is conical and regular. Suppose that $x_1+x_2=
y_1+y_2$ for some $x_i$, $y_j\in M$. If $x_1=0$, then there is a
refinement matrix
  \[
  \refmatr{y_1}{y_2}{x_1}00{x_2}{y_1}{y_2}
  \]
Similar refinements exist if $x_2$, $y_1$, or $y_2$ is zero.
Hence, we may assume that $x_i$, $y_j \in A$ for all $i$, $j$. In
the group $A$, we have $x_2= y_1+x_1^-+y_2$, and so
  \[
  \refmatr{y_1}{y_2}{x_1}{x_1}0{x_2}{y_1+x_1^-}{y_2}
  \]
is a refinement matrix.

(c$'$) Let $e\leq f$ in $\DI{M}$ and $x\in \AC{M}{e}$ such that
$x+f=f$. If $e=0$, then $x=0=e$. If $e\ne 0$, then $e=0_A$, whence
$f=0_A$ and $x\in A$. Since $A$ is a group, $x=0_A=e$ in this
case.

(d$'$) Let $e\leq f$ in $\DI{M}$ and $x\in \AC{M}{e}$, $y\in
\AC{M}{f}$ such that $x+f=my$ for some $m\in\NN$. If $e=0$, then
$x=0$, whence $x+f= f= m(e+f)$. If $e\ne 0$, then $e=0_A$, whence
$f=0_A$ and $x$, $y\in A$. In this case, $y\in \AC{M}{e}$, and
$x+f=m(y+f)$.
\end{proof}

\begin{definition}
We shall say that a \ram\ $M$ satisfies the \emph{embedding
condition\/}, abbreviated (emb), provided condition (c) of
Proposition~\ref{P:NecessLL} holds. Further, $M$ satisfies the
\emph{purity condition\/}, abbreviated (pur), provided $M$
satisfies condition (d) of the proposition.
\end{definition}

In view of the results above, any direct limit of finite direct
sums of monoids of the form $\znznz$ is a strongly periodic
conical refinement monoid satisfying (emb) and (pur). Our main
monoid-theoretic goal is to establish the converse statement
(Theorem~\ref{T:Main}).

We next investigate the structure of \ram s $M$, for which some
additional notation and terminology is helpful. Recall that
$a\propto b$ (for some $a$, $b\in M$) means that $a\leq mb$ for
some $m\in \NN$, and that $a\asymp b$ means that $a\propto
b\propto a$. Since $M$ is regular, $mb\leq b$ for all $m\in\NN$,
and so $a\propto b$ if and only if $a\leq b$. Thus, $a\asymp b$ if
and only if $a\leq b\leq a$. Similarly, $a\propto b$ if and only
if $\di{a}\propto\di{b}$, if and only if $\di{a}\leq \di{b}$.
Consequently, $a\asymp b$ if and only if $a$ and $b$ lie in the
same group $\AC{M}{e}$, for some idempotent $e\in \DI{M}$.

For any $a$, $b\in M$, the sum $\di{a}+\di{b}$ is an idempotent
with $\di{a}+\di{b} \asymp a+b$, whence $\di{a}+\di{b}=\di{a+b}$.
In particular, this shows that
$\AC{M}{e}+\AC{M}{f}\subseteq\AC{M}{e+f}$ for all idempotents $e$,
$f\in \DI{M}$. Now $e+f$ is the supremum of $e$ and $f$ in the
semilattice $\DI{M}$, but there need not exist an infimum. We do,
however, have a commutative diagram of abelian groups and group
homomorphisms as follows:
  \[
  \xymatrixcolsep{6pc}\xymatrixrowsep{1pc}
  \xymatrix{
   &\AC{M}{e} \ar[dr]^{y\mapsto y+f} \\
   {\varinjlim_{g\leq e,f}\AC{M}{g}}
   \ar[ur]^{x\mapsto x+e} \ar[dr]_{x\mapsto x+f} &&\AC{M}{e+f} \\
   &\AC{M}{f} \ar[ur]_{z\mapsto z+e}
  }
  \]
The resemblance of this diagram to a pullback behind a
Mayer-Vietoris sequence in homological algebra provides a
convenient name for the following monoid condition, which will be
our key to the refinement property in \ram s.

\begin{definition}
Let $M$ be a \ram. We shall say that $M$ satisfies the
\emph{Mayer-Vietoris property\/} (or \emph{MVP\/}, for short)
provided that, for all idempotents $e$, $f\in\DI{M}$,

(a) $\AC{M}{e}+\AC{M}{f}= \AC{M}{e+f}$.

(b) Whenever $u\in \AC{M}{e}$ and $v\in \AC{M}{f}$ with $u+f=
v+e$, there exists $w\in M$ such that $u=w+e$ and $v=w+f$. (Note
that necessarily $w\in \AC{M}{g}$ for some idempotent $g\leq e$,
$f$.)
\end{definition}

The following result is in some sense a version of Proposition~1
and Corollary~4 of~\cite{Dobb} with the finiteness assumption on
the monoid removed.

\begin{theorem}\label{T:CharRegRef}
A \ram\ $M$ is a refinement monoid if and only if $\DI{M}$ is a
distributive semilattice and $M$ satisfies the MVP.
\end{theorem}

\begin{proof}
$(\Longrightarrow)$: Suppose that $e_1+e_2= f_1+f_2$ for some
$e_i,f_j\in \DI{M}$. Refine this equation in $M$:
  \[
  \refmatr{f_1}{f_2}{e_1}{x_{11}}{x_{12}}{e_2}{x_{21}}{x_{22}}
  \]
Now if we set $g_{ij}= \di{x_{ij}}$ for all $i$, $j$, then $e_i=
x_{i1}+x_{i2} \in \AC{M}{g_{i1}+g_{i2}}$. Since $e_i$ is
idempotent, we obtain $g_{i1}+g_{i2} =e_i$ for $i=1$, $2$.
Similarly, $g_{1j}+g_{2j}=f_j$ for $j=1$, $2$, which shows that
$\DI{M}$ has refinement. Therefore $\DI{M}$ is a distributive
semilattice.

Now let $e$, $f\in \DI{M}$. We have already observed that
$\AC{M}{e}+\AC{M}{f}$ is contained in $\AC{M}{e+f}$. To prove the
reverse inclusion, consider an arbitrary element $a\in
\AC{M}{e+f}$. Note that $a= a+e+f$ and $a+a^-= e+f$. Take a
refinement of the second equation:
  \[
  \refmatr{e}{f}{a}{b}{c}{a^-}{u}{v}
  \]
Now $a= a+e+f= (b+e)+(c+f)$. Since $b+u=e$, we have $b\leq e$,
whence $b+e\asymp e$ and so $b+e\in \AC{M}{e}$. Similarly, $c+f\in
\AC{M}{f}$, and therefore $\AC{M}{e+f} \subseteq
\AC{M}{e}+\AC{M}{f}$. This establishes the first half of the MVP.

Given $u\in \AC{M}{e}$ and $v\in \AC{M}{f}$ with $u+f= v+e$, take
a refinement of this equation:
  \[
  \refmatr veuabfcd
  \]
Then $d\leq e$, $f$. Put $w:= a+d^-$. Then $w+e= a+d^-+b+d=
u+\di{d} =u$ because $\di{d}\leq e\leq u$, and $w+f= a+d^-+c+d=
v+\di{d}= v$ because $\di{d}\leq f\leq v$. Therefore $M$ satisfies
the MVP.

$(\Longleftarrow)$: Given $a_1+a_2= b_1+b_2$ in $M$, set
$e_i=\di{a_i}$ and $f_j= \di{b_j}$ for $i,j=1,2$, so that
$e_1+e_2= f_1+f_2$. Since $\DI{M}$ is distributive, it contains a
refinement
  \[
  \refmatr{f_1}{f_2}{e_1}{g_{11}}{g_{12}}{e_2}{g_{21}}{g_{22}}
  \]
By the MVP, each $\AC{M}{e_i}= \AC{M}{g_{i1}}+\AC{M}{g_{i2}}$, and
so each $a_i= c_{i1}+c_{i2}$ for some $c_{ij}\in \AC{M}{g_{ij}}$.
Note that $c_{1j}+c_{2j} \in \AC{M}{g_{1j}}+\AC{M}{g_{2j}}=
\AC{M}{f_j}$ for $j=1$, $2$, and that $(c_{11}+c_{21})+
(c_{12}+c_{22})= a_1+a_2= b_1+b_2$. Set $u:= c_{11}+c_{21} +b_1^-
\in \AC{M}{f_1}$ and $v:= b_2+c^-_{12}+c^-_{22} \in \AC{M}{f_2}$,
and observe that
  \[
  u+f_2= c_{11}+c_{21} +b_1^-+ c_{12}+ c^-_{12}+ c_{22}+ c^-_{22}=
  b_1+b_2+ b_1^-+ c^-_{12}+ c^-_{22}= v+f_1.
  \]
By the MVP, there exists an element $w\in M$ such that $u=w+f_1$
and $v=w+f_2$, and $w\in \AC{M}{h}$ for some idempotent $h\leq
f_1,f_2$. Then
  \begin{align*}
  c_{11}+c_{21}+w^- &= u+b_1+w^-= w+f_1+b_1+w^-= b_1+h+f_1 =b_1\\
  c_{12}+c_{22}+w &= c_{12}+c_{22}+ f_2+w= c_{12}+c_{22}+v=
  b_2+f_2 =b_2.
  \end{align*}

Since $h\leq f_1\leq e_1+e_2$, distributivity in $\DI{M}$ implies
that $h=h_1+h_2$ for some idempotents $h_i\leq e_i$. Applying the
MVP a final time, we obtain $w= w_1+w_2$ for some
$w_i\in\AC{M}{h_i}$. We check that
  \begin{align*}
  (c_{11}+w_1^-)+ (c_{21}+w_2^-) &= c_{11}+c_{21}+w^- =b_1\\
  (c_{12}+w_1)+ (c_{22}+w_2) &= c_{12}+c_{22}+w =b_2 \\
  (c_{i1}+w_i^-)+ (c_{i2}+w_i) &= a_i+h_i= a_i \qquad\qquad
  (i=1,\,2),
  \end{align*}
where the last equalities hold because $h_i\leq e_i\leq a_i$.
Therefore we have a refinement
  \begin{equation}
  \refmatr{b_1}{b_2}{a_1}{c_{11}+w_1^-}{c_{12}+w_1}
  {a_2}{c_{21}+w_2^-}{c_{22}+w_2}\tag*{\qed}
  \end{equation}
\renewcommand{\qed}{}
\medskip
\end{proof}

In particular, Theorem~\ref{T:CharRegRef} describes the conditions
needed to obtain refinement in a \ram\ $M(\Lambda,\cF)$
constructed from a semilattice $\Lambda$ and a functor~$\cF$ from
$\Lambda$ to abelian groups as in Section~\ref{S:Back}.  For
example, take $\Lambda= {\two}^2$, the Boolean monoid of subsets
of a $2$-element set. Viewed as a category obtained from a poset,
$\Lambda$ looks like this:
  \[
  \xymatrixcolsep{3pc}\xymatrixrowsep{0.25pc}
  
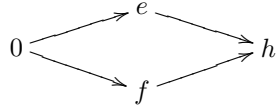
\xymatrix{
  &e \ar[dr] \\
  0 \ar[ur]\ar[dr] &&h \\
  &f \ar[ur]
  }
  \]
Suppose that $H$ is an abelian group with subgroups $E$, $F$, $G$
such that $G\subseteq E\cap F$. Then we can define a functor $\cF$
from $\Lambda$ to the category of abelian groups as follows:
  \[
  \xymatrixcolsep{3pc}\xymatrixrowsep{0.25pc}
  
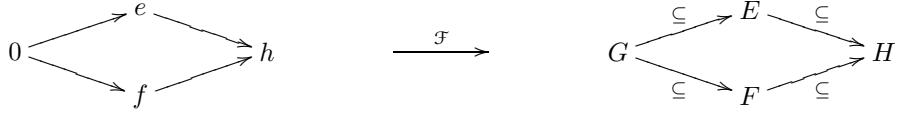
\xymatrix{
  &e \ar[dr] &&& &&E \ar[dr]^{\subseteq} \\
  0 \ar[ur]\ar[dr] &&h &\ar[r]^{\cF} &&G \ar[ur]^{\subseteq}
  \ar[dr]_{\subseteq} &&H \\
  &f \ar[ur] &&& &&F \ar[ur]_{\subseteq}
  }
  \]
Form the monoid $M= M(\Lambda,\cF)$. Then
Theorem~\ref{T:CharRegRef} says that $M$ has refinement if and
only if $E\cap F=G$ and $E+F=H$.

Because the group homomorphisms in the diagram above are
embeddings, the monoid $M$ is isomorphic to a submonoid of
$\Lambda\times H$, namely
  \[
  \bigl( \set{0} \times G \bigr) \sqcup \bigl( \set{e} \times E \bigr)
\sqcup \bigl( \set{f} \times F \bigr) \sqcup \bigl( \set{h} \times
H \bigr).
  \]
In fact, arbitrary \ram s with (emb) can be put into a similar
form, as follows.

\begin{theorem}\label{T:RAMstruct}
Let $M$ be a \ram\ satisfying the embedding condition. Then there
exist a semilattice $\Lambda$, an abelian group $G$, and subgroups
$G_e\subseteq G$ for all $e\in \Lambda$ such that
\begin{itemize}
\item[(a)] $G= \bigcup_{e\in\Lambda} G_e$.
\item[(b)] $G_e\subseteq G_f$ for all $e\leq f$ in $\Lambda$.
\item[(c)] $M$ is isomorphic to the submonoid
$\bigsqcup_{e\in\Lambda} \bigl( \{e\}\times G_e \bigr) \subseteq
\Lambda\times G$.
\end{itemize}
The monoid $M$ is a refinement monoid if and only if
\begin{itemize}
\item[(a$'$)] $\Lambda$ is distributive.
\item[(b$'$)] $G_e+G_f= G_{e+f}$ for all $e$, $f\in \Lambda$.
\item[(c$'$)] $G_e\cap G_f= \bigcup_{g\in\Lambda,\, g\leq e,f} G_g$
for all $e$, $f\in \Lambda$.
\end{itemize}
Moreover, $M$ is conical if and only if
\begin{itemize}
\item[(d$'$)] $G_0=\{0\}$,
\end{itemize}
and $M$ satisfies the purity condition if and only if
\begin{itemize}
\item[(e$'$)] $G_e$ is a pure subgroup of $G$ for all $e\in\Lambda$.
\end{itemize}
\end{theorem}

\begin{proof}
Set $\Lambda= \DI{M}$, and for $e\leq f$ in $\Lambda$, let
$\phi_{e,f}\colon\AC{M}{e} \rightarrow \AC{M}{f}$ denote the
homomorphism $x\mapsto x+f$. The collection of groups $\AC{M}{e}$
and transition maps $\phi_{e,f}$ forms a direct system in the
category of abelian groups. Let $G$ be the direct limit of this
system, with limiting maps $\eta_e \colon \AC{M}{e} \rightarrow G$
for $e\in\Lambda$, and set $G_e= \eta_e(\AC{M}{e})$ for $e\in
\Lambda$. Conditions (a) and (b) are clear, and the isomorphism
required in (c) is given by the rule $a\mapsto
(\di{a},\eta_{\di{a}}(a))$.

It follows from Theorem~\ref{T:CharRegRef} that $M$ is a
refinement monoid if and only if (a$'$), (b$'$), (c$'$) hold, and
the remaining equivalences are clear. (Note that (e$'$) is
equivalent to the statement that $G_e$ is pure in $G_f$ whenever
$e\le f$ in $\Lambda$.)
\end{proof}

For certain applications, it is useful to be able to restrict to
strongly periodic monoids in which the orders of the elements are
controlled, as follows.

Recall that a \emph{generalized integer\/} or \emph{supernatural
number\/} is a formal product of nonnegative powers of the
positive prime integers, thus
  \[
  \prod_p p^{t(p)} = 2^{t(2)} 3^{t(3)} 5^{t(5)} \cdots p^{t(p)}
  \cdots\,,
  \]
where each exponent $t(p)\in \{0\}\cup\NN\cup\{\infty\}$. If $\fm=
\prod_p p^{s(p)}$ and $\fn= \prod_p  p^{t(p)}$ are generalized
integers, the statement $\fm\mid\fn$ means that $s(p)\leq t(p)$
for all primes $p$. Ordinary positive integers are treated as
generalized integers in the obvious manner.

\begin{definition} For any \ram\ $M$ and generalized
integer $\fm$, we set
  \[
  M[\fm]= \setm{x\in M}{(m+1)x=x \text{\ for some positive integer\ }
  m\mid\fm}.
  \]

Note that $M[\fm]$ is a submonoid of $M$ containing $\DI{M}$, and
that it is also a semilattice of groups, since the sets
  \begin{align*}
  M[\fm]\cap \AC{M}{e} &= \setm{x\in M}{mx=e \text{\ for some
positive integer\ } m\mid\fm}
  \end{align*}
are subgroups of $\AC{M}{e}$ for each $e\in \DI{M}$.
\end{definition}

\begin{proposition}\label{P:M[m]transf}
Let $M$ be a regular refinement monoid satisfying the embedding
and purity conditions, and let $\fm$ be a generalized integer.
Then $M[\fm]$ is a regular refinement monoid satisfying the
embedding and purity conditions.
\end{proposition}

\begin{proof} We have already observed that $M[\fm]$ is a
semilattice of groups, and that $\DI{M[\fm]}= \DI{M}$, whence
$\DI{M[\fm]}$ is a distributive semilattice. It is clear that
(emb) passes from $M$ to $M[\fm]$.

Let $e$, $f$, $g$ be idempotents in $M$ with $e+f=g$. If $z\in
\AC{M[\fm]}{g}$, then $mz=g$ for some positive integer $m\mid\fm$.
By the MVP, $z=b+c$ for some $b\in \AC{M}{e}$ and $c\in
\AC{M}{f}$. Add $mc^-$ to both sides of the equation $mb+mc= g$,
to obtain $mb+f= mc^-+e$. The MVP now implies that there exists an
element $w\in M$ such that $mb=w+e$ and $mc^- = w+f$; moreover,
$w\in \AC{M}{h}$ for some idempotent $h\leq e$, $f$. Since
$w+e=mb$, it follows from (pur) and (emb) that $w=mv$ for some
$v\in \AC{M}{h}$. Set $v'= (v+e)^- \in \AC{M}{e}$. Since $mb= w+e=
m(v+e)$, the element $b+v'\in \AC{M}{e}$ satisfies $m(b+v')=e$.
Moreover, $c+v+f\in \AC{M}{f}$, and $mc^- = w+f= mv+f$ implies
$m(c+v+f)=f$. Finally,
  \begin{align*}
  (b+v')+(c+v+f) &= b+c+(v+e)^-+v+f\\
  &=z+(v+e)^-+(v+e)+f=z+e+f= z.
  \end{align*}
Thus, $\AC{M[\fm]}{g}=\AC{M[\fm]}{e}+ \AC{M[\fm]}{f}$. Now suppose
that $u\in\AC{M[\fm]}{e}$ and $v\in \AC{M[\fm]}{f}$ with
$u+f=v+e$. By the MVP in $M$, there exists an element $w\in M$
such that $u=w+e$ and $v=w+f$. Put $h=\di{w}$, and choose
$m\in\NN$, with $m\mid\fm$, such that $mu=e$ and $mv=f$. Since
$mw+e=mu=e$, (emb) implies that $mw=h$, so that $w\in
\AC{M[\fm]}{h}$. This shows that $M[\fm]$ satisfies the MVP.
Therefore, by Theorem~\ref{T:CharRegRef}, $M[\fm]$ is a refinement
monoid.

Let $e\leq f$ be idempotents in $M$, and consider elements $x\in
\AC{M[\fm]}{e}$ and $y\in \AC{M[\fm]}{f}$ such that $x+f= ny$ for
some $n\in\NN$. Choose $m\in\NN$, with $m\mid\fm$, such that
$mx=e$ and $my=f$, and let $d= \text{GCD}(m,n)$. Then $m=m'd$ and
$n=n'd$ for some $m'$, $n'\in \NN$, and $\text{GCD}(m',n')=1$.
Note that $m'x+f= m'ny= n'my= f$, whence $m'x=e$ by (emb). Now
$x+f= d(n'y)$ with $n'y\in \AC{M}{f}$. Using (pur) and (emb) in
$M$, we obtain an element $z\in \AC{M}{e}$ such that $x=dz$.
Moreover, $mz= m'x= e$, and so $z\in M[\fm]$. Since $n'$ and $m'$
are relatively prime, there exists $n^*\in \NN$ such that $n^*n'
\equiv 1 \pmod{m'}$, whence $n^*n \equiv d \pmod{m}$, and so
$n^*nz=dz$. Thus $x= dz= n(n^*z)$ with $n^*z\in \AC{M[\fm]}{e}$,
which establishes (pur) in~$M[\fm]$.
\end{proof}

\section{Direct limits}\label{S:DirLim}

Since our aim is to express certain monoids as direct limits of
appropriate building blocks, it is helpful to set down general
conditions for such direct limits at the outset. We shall use the
following version of \cite[Lemma~3.4]{GW}, which many readers will
recognize as an analogue of a key step in other classification
results. It is a monoid-theoretical version of Shannon's result
\cite[Theorem~2]{Shan}. For a map $\phi\colon X\to Y$, we put
  \[
  \ker\phi=\setm{(x,y)\in X\times X}{\phi(x)=\phi(y)}.
  \]

\begin{lemma}\label{L:TrLemma}
Let $\cB$ be a class of finite abelian monoids which is closed
under finite direct sums and let $M$ be an abelian monoid. Then
$M$ is a direct limit of monoids from $\cB$ if and only if the
following two conditions are satisfied:
  \begin{itemize}
  \item[(1)] For each $x\in M$, there exist $B\in\cB$ and a
homomorphism $\phi\colon B \rightarrow M$ such that $x\in\phi(B)$.
  \item[(2)] For any $B\in\cB$ and any homomorphism
$\phi\colon B\rightarrow M$, there exist $B'\in \cB$ and
homomorphisms $\xymatrix{B\ar[r]^-{\psi} & B'\ar[r]^-{\phi'} & M}$
such that $\phi'\psi= \phi$ and $\ker\phi=\ker\psi$.
  \end{itemize}
\end{lemma}

\begin{proof}
The given conditions clearly imply the two hypotheses of
\cite[Lemma~3.4]{GW}, hence they imply that $M$ is a direct limit
of members of $\cB$.

Conversely, suppose that $M=\varinjlim_{i\in I}B_i$, a direct
limit with all $B_i$ in $\cB$, transition maps $f_{ij}\colon
B_i\to B_j$, and limiting maps $f_i\colon B_i\to M$, for all
$i\leq j$ in the directed partially ordered set $I$. As
$M=\bigcup_{i\in I} f_i(B_i)$, Condition~(1) is satisfied. Now let
$\phi\colon B\to\nobreak M$ be a monoid homomorphism, with
$B\in\cB$. Since $B$ is finite, $\phi(B) \subseteq f_i(B_i)$ for
some $i\in I$. Choose elements $x_b \in B_i$ such that $f_i(x_b)=
\phi(b)$ for all $b\in B$ and $x_0=0$. For all $c,d\in B$, we have
$f_i(x_c+x_d)= \phi(c+d)= f_i(x_{c+d})$. By finiteness, there is
some $j\in I$, with $j\geq i$, such that $f_{ij}(x_c+x_d)=
f_{ij}(x_{c+d})$ for all $c,d\in B$. Now replace $i$ by $j$ and
each $x_b$ by $f_{ij}(x_b)$. This allows us to assume, without
loss of generality, that $x_c+x_d= x_{c+d}$ for all $c,d\in B$.
Hence, there is a monoid homomorphism $\psi\colon B\rightarrow
B_i$, given by $\psi(b)= x_b$, such that $f_i\psi= \phi$. For each
$(x,y)\in \ker\phi$, we have $f_i\psi(x)= f_i\psi(y)$, and so
there is some $k\in I$, with $k\geq i$, such that $f_{ik}\psi(x)=
f_{ik}\psi(y)$ for all $(x,y)\in \ker\phi$. Now replace $i$
and~$\psi$ by~$k$ and $f_{ik}\psi$. This allows to assume that
$\ker\phi \subseteq \ker\psi$. Since the reverse inclusion follows
from $f_i\psi= \phi$, we conclude that (2) above is satisfied with
$B'= B_i$ and $\phi'= f_i$.
\end{proof}

In an arbitrary category admitting all direct limits (in
categorical language, directed colimits), the class of all direct
limits of members from a given class is not necessarily closed
under direct limits -- even in case the category we are starting
with is a partially ordered set! However, strengthening the
assumptions leads to the following useful positive result.

\begin{corollary}\label{C:TrCor}
Let $\cB$ be a class of finite abelian monoids which is closed
under finite direct sums. Then the class of all direct limits of
monoids from $\cB$ is closed under direct limits.
\end{corollary}

\begin{proof}
Denote by $\cL$ the class of all direct limits of monoids from
$\cB$. Let $M=\varinjlim_{i\in I}M_i$, a direct limit with all
$M_i\in\cL$, transition maps $f_{ij}\colon M_i\to M_j$ and
limiting maps $f_i\colon M_i\to M$, for all $i\leq j$ in the
directed partially ordered set $I$. Since the $M_i$ satisfy
Condition~(1) of Lemma~\ref{L:TrLemma} and $M=\bigcup_{i\in I}
f_i(M_i)$, we see that~$M$ satisfies Condition~(1) of
Lemma~\ref{L:TrLemma}. Now let $\phi\colon B\to M$ be a monoid
homomorphism, with $B\in\cB$. Since $B$ is finite, we see as in
the proof of Lemma~\ref{L:TrLemma} that there are $i\in I$ and a
monoid homomorphism $\phi'\colon B\to M_i$ such that
$\phi=f_i\phi'$ and $\ker\phi=\ker\phi'$. Since $M_i\in\cL$,
Lemma~\ref{L:TrLemma} shows that there exists $B'\in\cB$ together
with monoid homomorphisms $\psi\colon B\to B'$ and $\phi''\colon
B'\to M_i$ such that $\phi'=\phi''\psi$ and $\ker\phi'=\ker\psi$.
Therefore, $\phi= (f_i\phi'')\psi$ with $f_i\phi''\colon B'\to M$
and $\ker\phi=\ker\psi$. Using Lemma~\ref{L:TrLemma} again, we
conclude that $M$ belongs to $\cL$.
\end{proof}

\begin{remark}\label{Rk:FinPres}
Both Lemma~\ref{L:TrLemma} and Corollary~\ref{C:TrCor} can be
extended to the case where all members of $\cB$ are \emph{finitely
generated\/} monoids. To obtain this, we observe that in the proof
of Lemma~\ref{L:TrLemma}, the monoid $B/{\ker\phi}$ is finitely
generated, thus, by Redei's Theorem (see \cite{Redei}, or
\cite{Frey} for a simple proof), finitely presented.
\end{remark}

For the remainder of the paper, we restrict $\cB$ to be the class
of finite direct sums of monoids of the form $\znznz$ for
$n\in\NN$, and we let $\cL$ denote the class of all direct limits
of monoids from $\cB$. Further,  write $\Rep$ for the class of all
strongly periodic conical refinement monoids satisfying the
conditions (emb) and (pur). It follows from
Proposition~\ref{P:NecessLL} that $\cL$ is contained in $\Rep$,
and the main goal of Sections~\ref{S:FinMon} and~\ref{S:FinRep} is
to prove the reverse inclusion.

\begin{lemma}\label{L:DirProd}
The class $\cL$ is closed under direct limits, finite direct sums,
and retracts.
\end{lemma}

\begin{proof} Corollary~\ref{C:TrCor} implies that $\cL$ is closed
under direct limits, and it is straightforward to verify that
$\cL$ is closed under finite direct sums.

Now consider a monoid $M$ which is a retract of a monoid
$M'\in\cL$, that is, there are morphisms $\eps\colon M\to M'$ and
$\mu\colon M'\to M$ such that $\mu\eps=\id_M$. Put $\rho=\eps\mu$,
and observe that $\rho^2=\rho$ and $\mu\rho=\mu$. We claim that
$M$ is the direct limit of the sequence
  \[
  {
  \def\labelstyle{\displaystyle}
  \xymatrix{
  M'\ar[r]^-{\rho} & M'\ar[r]^-{\rho} & M'\ar[r]^-{\rho} &
  \ar@{}|{\displaystyle\cdots}[r] & ,
  }
  }
  \]
with constant limiting morphism $\mu\colon M'\to M$. Suppose that
we have a monoid $C$ and morphisms $\varphi_n\colon M'\rightarrow
C$ for $n\in\NN$ such that $\varphi_n=\varphi_{n+1}\rho$ for all
$n$. Since $\rho$ is idempotent, $\varphi_n=\varphi_0$ for all
$n$, and so $\varphi_0\eps$ is the unique morphism $\psi\colon
M\to C$ such that $\psi\mu=\varphi_0$. This establishes the claim,
and since $\cL$ is closed under direct limits, we conclude that
$M\in\cL$.
\end{proof}

\begin{corollary}\label{C:nzAincL}
For any finite abelian group $A$, the monoid $\nz{A}$ belongs to
$\cL$.
\end{corollary}

\begin{proof}
By the fundamental structure theorem of finite abelian groups,
$A=\bigoplus_{i=1}^nA_i$ for some finite cyclic groups $A_i$. Now
set $M=\bigoplus_{i=1}^n A_i^{\sqcup0}$, and note that the
inclusion map $A\into M$ extends to a unique monoid embedding
$\eps\colon\nz{A}\into M$.

For $i=1,\dots,n$, the canonical injection $A_i\into A$ extends to
a unique monoid embedding $\mu_i\colon A_i^{\sqcup0} \into\nz{A}$.
The maps $\mu_i$ induce a monoid homomorphism $\mu\colon
M\to\nz{A}$ given by the rule
$\mu(a_1,\dots,a_n)=\sum_{i=1}^n\mu_i(a_i)$. It is clear that
$\mu\eps$ is the identity map on $\nz{A}$, whence $\nz{A}$ is a
retract of $M$. Therefore, by Lemma~\ref{L:DirProd},
$\nz{A}\in\cL$.
\end{proof}

\section{Finite monoids}\label{S:FinMon}

The first major step towards our main result is to show that every
finite monoid from $\Rep$ belongs to $\cL$. We do this in the
present section, after recalling some facts about \jirr\ elements
in semilattices.

Every finite semilattice is, of course, a lattice, and it is
distributive as a semilattice if and only if it is distributive as
a lattice. A nonzero (i.e., non-minimum) element $p$ in a
semilattice $S$ is \emph{\jirr\/} if $p$ is not the supremum of
any pair of elements less than $p$, that is, if $p=x\vee y$
implies that $p\in\set{x,y}$, for any $x$, $y\in S$. We denote by
$\J(S)$ the set of all \jirr\ elements of $S$, and, for each $a\in
S$, we put $\J_S(a)=\setm{p\in\J(S)}{p\leq a}$. It is well-known
(see \cite[Exercise~I.6.13]{GLT2}) that in case $S$ is finite,
every element of $S$ is the supremum of the \jirr\ elements it
dominates, that is, $a=\bigvee\J_S(a)$ for all $a\in S$.
Furthermore, an element $p\in S$ is \jirr\ if and only if $p$ has
a unique \emph{lower cover\/}, that is, an element $x<p$ in $S$
such that no $y\in S$ satisfies $x<y<p$. In that case we shall
denote by~$p_*$ the unique lower cover of~$p$.

The following lemma is folklore.

\begin{lemma}\label{L:pdagger}
For every \jirr\ element $p$ in a finite distributive lattice $D$,
there exists a unique largest $u\in D$ such that $p\nleq u$.
\end{lemma}

\begin{proof}
Since $D$ is distributive and $p$ is \jirr, $p\nleq x$ and $p\nleq
y$ implies that $p\nleq x\vee y$, for any $x$, $y\in D$. Set $u=
\bigvee \setm{x\in D}{p\nleq x}$.
\end{proof}

The element $u$ of Lemma~\ref{L:pdagger} is traditionally denoted
by $p^{\dagger}$.

For an abelian group $G$, let us denote by $\Sub G$ the lattice of
all subgroups of~$G$. The following lemma is also folklore. It is
valid in the much more general context of a homomorphism from a
finite distributive lattice to a modular lattice with zero.

\begin{lemma}\label{L:DistrGL}
Let $G$ be an abelian group, $D$ a finite distributive lattice,
$f\colon D\to\Sub G$ a lattice homomorphism, and
$\famm{H_p}{p\in\J(D)}$ a family of subgroups of $G$ such that
$f(p)= f(p_*)\oplus H_p$ for all $p\in\J(D)$. Then
  \[
  f(a)=f(0)\oplus\bigoplus_{p\in\J_D(a)} H_p
  \]
for all $a\in D$.
\end{lemma}

\begin{proof}
We argue by induction on $a$. As the result is trivial for $a=0$
(in which case $\J_D(a)$ is empty), we only deal with the
induction step. Let $b$ be a lower cover of $a$ in $D$ and let
$p\leq a$ be minimal with respect to the property $p\nleq b$. Then
$p$ is \jirr, and, by the minimality statement, $p_*\leq b$.
Hence, $p\wedge b=p_*$ and $p\vee b=a$. For any $q\in\J(D)$ such
that $q\leq a$, it follows from the \jirry\ of $q$ and the
distributivity of $D$ that either $q\leq b$ or $q\leq p$. If
$q\nleq b$, then $q\leq p$, and $q<p$ is ruled out because that
would imply $q\leq p_*\leq b$, a contradiction. Hence, we have
proved the statement
  \begin{equation}\label{Eq:JD(a)}
  \J_D(a)=\J_D(b)\cup\set{p}.
  \end{equation}
Now we compute:
  \[
  f(b)+H_p=f(b)+f(p_*)+H_p=f(b)+f(p)=f(b\vee p)=f(a)
  \]
because $f(p_*) \subseteq f(b)$, while
  \[
  f(b)\cap H_p=f(b)\cap f(p)\cap H_p=f(b\wedge p)\cap H_p=
  f(p_*)\cap H_p=\set{0}
  \]
because $H_p \subseteq f(p)$. Therefore, $f(a)=f(b)\oplus H_p$,
and thus, by \eqref{Eq:JD(a)} and the induction hypothesis,
$f(a)=f(0)\oplus \bigoplus_{q\in\J_D(a)} H_q$.
\end{proof}

\begin{proposition}\label{P:FinRep}
Any finite monoid in $\Rep$ belongs to $\cL$.
\end{proposition}

\begin{proof}
Let $M$ be a finite monoid in $\Rep$. In view of
Theorem~\ref{T:RAMstruct}, we may assume that
  \[
  M= \bigsqcup_{e\in\Lambda}
\bigl( \{e\}\times G_e \bigr) \subseteq \Lambda\times G
  \]
for some finite semilattice $\Lambda$ and some finite abelian
group $G$ with subgroups $G_e$ (for $e\in\Lambda$) satisfying the
conditions (a), (b), and (a$'$)--(e$'$) of the theorem. Finally,
since~$\Lambda$ is finite, it is a distributive lattice, and
condition (c$'$) implies that $G_e\cap G_f= G_{e\wedge f}$ for all
$e$, $f\in\Lambda$. Note that the rule $e\mapsto G_e$ provides a
lattice homomorphism $\Lambda \rightarrow \Sub G$.

For any $p\in\J(\Lambda)$, the group $G_{p_*}$ is a finite, pure
subgroup of $G_p$, and so, by Kulikov's Theorem (see
\cite[Theorem~27.5]{FuchsI}), $G_p=G_{p_*}\oplus H_p$ for some
subgroup $H_p$ of $G_p$. Lemma~\ref{L:DistrGL} thus yields that
  \begin{equation}\label{Eq:GaDirSum}
  G_e=\bigoplus_{p\in\J_\Lambda(e)} H_p
  \end{equation}
for all $e\in\Lambda$. In particular, taking $e=1$ (the maximum
element of $\Lambda$), we obtain $G= \bigoplus_{p\in\J(\Lambda)}
H_p$. Let $\pi_q \colon G\rightarrow H_q$, for $q\in\J(\Lambda)$,
denote the projections corresponding to this direct sum.

We next define maps $\eps_p\colon M\to\nz{G}$ and
$\mu_p\colon\nz{G}\to M$, for $p\in\J(\Lambda)$, by the rules
  \begin{align*}
  \eps_p(e,x)&=\begin{cases}
  \pi_p(x)&(p\leq e)\\
  \zo&(p\nleq e)
  \end{cases}&
  \mu_p(y)&=\begin{cases}
  (0,0)&(y=\zo)\\
  (p,\pi_p(y))&(y\in G).
  \end{cases}
  \end{align*}
It is clear that $\mu_p$ is a monoid homomorphism, and we claim
that $\eps_p$ is one as well. Hence, we must show that
  \begin{equation}\label{Eq:epsphom}
  \eps_p(e\vee f,\, x+y) = \eps_p(e,x)+ \eps_p(f,y)
  \end{equation}
for all $(e,x),(f,y) \in M$. If $p\leq e$ and $p\leq f$, then both
sides of \eqref{Eq:epsphom} equal $\pi_p(x+y)$, while if $p\nleq
e$ and $p\nleq f$, both sides are zero. If $p\nleq e$ but $p\leq
f$, then in view of~\eqref{Eq:GaDirSum}, $\pi_p(x)= 0$ (because
$p\notin \J_\Lambda(e)$), whence both sides of \eqref{Eq:epsphom}
equal $\pi_p(y)$. A symmetric observation covers the remaining
situation, and thus \eqref{Eq:epsphom} holds in all cases.

Finally, we define homomorphisms $\eps\colon
M\to(\nz{G})^{\J(\Lambda)}$ and
$\mu\colon(\nz{G})^{\J(\Lambda)}\to M$ by the rules
  \begin{align*}
  \eps(e,x)&= \bigl( \eps_p(e,x) \bigr)_{p\in\J(\Lambda)}
  &\mu \bigl( (y_p)_{p\in\J(\Lambda)} \bigr) &= \sum_{p\in\J(\Lambda)}
  \mu_p(y_p).
  \end{align*}
For any nonzero $(e,x)\in M$, we compute that
  \begin{align*}
  \mu\eps(e,x) &= \sum_{p\in\J(\Lambda)} \mu_p \eps_p(e,x)=
\sum_{p\in\J_\Lambda(e)} \mu_p \pi_p(x) \\
  &= \sum_{p\in\J_\Lambda(e)} (p,\pi_p(x))= \Bigl( e,
\sum_{p\in\J_\Lambda(e)} \pi_p(x) \Bigr)= (e,x),
  \end{align*}
where the final equality comes from \eqref{Eq:GaDirSum}. Thus,
$\mu\eps=\id_M$, and so $M$ is a retract of
$(\nz{G})^{\J(\Lambda)}$. We conclude from
Corollary~\ref{C:nzAincL} and Lemma~\ref{L:DirProd} that
$M\in\cL$.
\end{proof}

\begin{remark}
The direct limits that exist by virtue of
Proposition~\ref{P:FinRep} necessarily involve systems of
non-injective homomorphisms, even in the case of semilattices --
while every distributive semilattice is a direct limit of finite
Boolean semilattices \cite[Theorem 6.6]{GW}, most distributive
semilattices are not directed unions of finite Boolean
subsemilattices. This is just because finite distributive
semilattices need not be Boolean, the three-element chain
$\set{0,1,2}$ being the simplest example. This semilattice can be
expressed as a direct limit of copies of $\two^2$; see
\cite[Example 6.8]{GW}.
\end{remark}

\section{Characterization of the monoids in $\Rep$}
\label{S:FinRep}

Because of Proposition~\ref{P:FinRep}, we will be able to conclude
that $\Rep=\cL$ once we show that every monoid in $\Rep$ is a
direct limit of finite members of $\Rep$. In fact, we will show
that monoids in $\Rep$ are directed unions of finite submonoids
from~$\Rep$. This also provides a generalization of Pudl\'ak's
result, \cite[Fact~4, p.~100]{Pudl85}, that every distributive
semilattice is the directed union of its finite distributive
subsemilattices.

\begin{theorem}\label{T:FinRep}
Each monoid $M$ in $\Rep$ is the directed union of those finite
submonoids of $M$ which belong to $\Rep$.
\end{theorem}

\begin{proof}
We must show that any finite subset $X$ of $M$ is contained in
some finite submonoid of $M$ lying in $\Rep$. For convenience,
assume that $0\in X$. We first reduce to the case where there is a
bound on the orders of the elements of $M$, by observing that $M$
is the directed union of all $M[m]$, for $m\in\NN$; thus,
$X\subseteq M[m]$ for some $m$. By Proposition~\ref{P:M[m]transf},
$M[m]\in\Rep$, and so we may replace $M$ by $M[m]$.

Hence, we may assume that $(m+1)x=x$ for all $x\in M$, where $m$
is a fixed positive integer. We start as in the proof of
Proposition~\ref{P:FinRep}. By Theorem~\ref{T:RAMstruct}, we may
assume that
  \[
  M= \bigsqcup_{e\in\Lambda}
\bigl( \{e\}\times G_e \bigr) \subseteq \Lambda\times G
  \]
for some distributive semilattice $\Lambda$ and some abelian group
$G$ with subgroups $G_e$ satisfying all the conditions of the
theorem.

Next, we set $G_A=\bigcup_{e\in A} G_e$ for every ideal $A$ of
$\Lambda$. Observe that the union defining $G_A$ is directed, and
that $G_{[0,e]}=G_e$ for all $e\in \Lambda$. Hence, if $A\subseteq
B$ in $\Id \Lambda$, then $G_A$ is a pure subgroup of $G_B$. Since
$mG_A=\set{0}$, it follows from Kulikov's Theorem that $G_A$ must
be a direct summand of $G_B$. Notice also that $G_A+G_B= G_{A\vee
B}$ and $G_A\cap G_B= G_{A\cap B}$ for arbitrary $A$, $B\in \Id
\Lambda$. Thus, the rule $A\mapsto G_A$ defines a lattice
homomorphism $\Id\Lambda \rightarrow \Sub G$.

Write the elements $x\in X$ in the form $x=(e_x,g_x)\in M$. Denote
by $\bD$ the sublattice of $\Id \Lambda$ generated by the
principal ideals $[0,e_x]$ for $x\in X$. Since $\Id \Lambda$ is
distributive, $\bD$ is finite (in fact, $|\bD|\leq2^{2^{|X|}}$).
Moreover, the ideal $\set{0}$ belongs to~$\bD$ because $0\in X$.
For each $P\in\J(\bD)$, choose a subgroup $H_{P}$ of $G_P$ such
that $G_P= G_{P_*}\oplus H_{P}$, where $P_*$ denotes the unique
lower cover of $P$ in the lattice~$\bD$. Lemma~\ref{L:DistrGL} now
implies that
  \[
  G_A=\bigoplus_{P\in\J_{\bD}(A)} H_{P}
  \]
for all $A\in\bD$. In particular, taking $A$ to be the largest
element, say $I$, of $\bD$, we obtain $G_I=
\bigoplus_{P\in\J(\bD)} H_{P}$.

For each $x\in X$, we have
  \[
  g_x\in G_{e_x}= G_{[0,e_x]}= \bigoplus_{P\in\J_{\bD}([0,e_x])}H_{P}.
  \]
Since $X$ is finite, there exist finitely generated subgroups
$H'_{P}\subseteq H_{P}$ for $P\in\J(\bD)$ such that
  \begin{equation}\label{Eq:H'fplarge}
  g_x\in\bigoplus_{P\in\J_{\bD}([0,e_x])} H'_{P}
  \end{equation}
for $x\in X$. Since each $mH_P=0$, the groups $H'_P$ are all
finite. Define finite subgroups
  \begin{equation}\label{Eq:DefG'fa}
  G'_A=\bigoplus_{P\in\J_{\bD}(A)} H'_{P} \subseteq G_A
  \end{equation}
for all $A\in\bD$. Observe that
  \begin{equation}\label{Eq:G'pluscap}
  G'_A+G'_B = G'_{A+B} \quad\text{and}\quad G'_A\cap G'_B = G'_{A\cap
B} \quad\text{for all\ } A,\,B\in \bD,
  \end{equation}
and that
  \begin{equation}\label{Eq:G'fapure}
  G'_A \text{\ is a pure subgroup of\ } G'_B \quad\text{for all\ }
A\subseteq B \text{\ in\ } \bD.
  \end{equation}

For each $x\in X$, since $[0,e_x]$ is the supremum of all \jirr\
elements of $\bD$ below it, there are elements $u^x_{P}\in P$, for
$P\in\J_{\bD}([0,e_x])$, such that $e_x=
\bigvee_{P\in\J_{\bD}([0,e_x])} u^x_{P}$. Setting
$u_{P}=\bigvee_{x\in X,\,[0,e_x]\supseteq P} u^x_{P}$ for
$P\in\J(\bD)$, we obtain that $u_{P}\in P$ and
  \begin{equation}\label{Eq:Bndaxfp}
  e_x= \bigvee_{P\in\J_{\bD}([0,e_x])}u_{P}
  \end{equation}
for all $x\in X$. Since each $G'_{P}$ is a finite subset of the
directed union $G_{P}=\bigcup_{e\in P} G_e$, there exist elements
$v_{P}\in P$ such that $G'_{P}\subseteq G_{v_{P}}$ for all
$P\in\J(\bD)$. Finally, for each $P\in\J(\bD)$, recall the
notation $P^\dagger$ for the unique largest element of $\bD$ not
containing $P$ (see Lemma~\ref{L:pdagger}), choose $w_{P}\in
P\setminus P^{\dagger}$, and put $\psi(P)=u_{P}\vee v_{P}\vee
w_{P}$. We define a map $\varphi\colon\bD\to \Lambda$ by the rule
  \[
  \varphi(A)=\bigvee_{P\in\J_{\bD}(A)} \psi(P),
  \]
and we claim that
\begin{itemize}
\item[(1)] $\varphi$ is a semilattice embedding.
\item[(2)] $\varphi(\bD)$ is a finite distributive subsemilattice of
$\Lambda$.
\item[(3)] $\varphi(A)\in A$ for all $A\in\bD$.
\item[(4)] $\varphi([0,e_x])=e_x$ for all $x\in X$.
\end{itemize}
The third statement is clear since $\psi(P)\in P$ for all
$P\in\J(\bD)$. In particular, $\varphi(\set{0})= 0$. It is also
clear that $\varphi$ is a semilattice homomorphism. To finish the
proof of (1), consider $A$, $B\in\bD$ such that $A \not\subseteq
B$. There exists $P\in\J(\bD)$ such that $P\subseteq A$ but
$P\not\subseteq B$, and then $B\subseteq P^{\dagger}$. From
$P\subseteq A$ it follows that $w_{P}\leq\varphi(A)$. On the other
hand, from $w_{P}\notin P^{\dagger}$ it follows that $w_{P}\notin
B$, and so $w_{P}\nleq\varphi(B)$. Therefore,
$\varphi(A)\nleq\varphi(B)$, and (1) is proved. It now follows
that $\varphi(\bD)$ is a finite subsemilattice of $\Lambda$,
isomorphic to $\bD$ and hence distributive, establishing~(2).
Finally, for $x\in X$, it follows from (3) that
$\varphi([0,e_x])\leq e_x$. On the other hand,
  \[
  \varphi([0,e_x]) =
  \bigvee_{P\in\J_{\bD}([0,e_x])} \psi(P)
  \geq \bigvee_{P\in\J_{\bD}([0,e_x])} u_{P} = e_x
  \]
by \eqref{Eq:Bndaxfp}, and (4) is proved.

Now we set $N=\bigsqcup_{A\in\bD} \bigl( \set{\varphi(A)}\times
G'_{A} \bigr) \subseteq \Lambda\times G$. In  view of
\eqref{Eq:G'pluscap}, $N$ is a finite submonoid of $\Lambda\times
G$. Since
  \[
  G'_{A} =\sum_{P\in\J_{\bD}(A)} H'_{P}
  \subseteq \sum_{P\in\J_{\bD}(A)} G_{v_{P}}
  \subseteq \sum_{P\in\J_{\bD}(A)} G_{\psi(P)} =
  G_{\varphi(A)}
  \]
for all $A\in\bD$, we see that $N \subseteq M$. By (2),
$\DI{N}\cong \varphi(\bD)$ is a (finite) distributive semilattice.
It now follows from \eqref{Eq:G'pluscap} and
Theorem~\ref{T:RAMstruct} that $N$ is a refinement monoid. It is
clear that $N$ is conical and satisfies (emb), and $N$ satisfies
(pur) by~\eqref{Eq:G'fapure}. Thus, $N$ belongs to $\Rep$.

Finally, for every $x\in X$,
  \[
  g_x\in\bigoplus_{P\in\J_{\bD}([0,e_x])} H'_{P} =
  G'_{[0,e_x]}
  \]
by \eqref{Eq:H'fplarge} and \eqref{Eq:DefG'fa}, whence
$x=(e_x,g_x)\in N$. Therefore, $X$ is contained in $N$.
\end{proof}

\begin{remark}\label{Rk:dontPudlak}
It is tempting to try to reduce the proof of
Theorem~\ref{T:FinRep} to the case where $\Lambda$ is finite, by
applying Pudl\'ak's result. After putting $M$ into the form given
by Theorem~\ref{T:RAMstruct}, we can choose a finite set
$E\subseteq\Lambda$ such that $X\subseteq \bigsqcup_{e\in E}
\bigl(\set{e}\times G_e\bigr)$; then, by Pudl\'ak's result,
$\Lambda$ has a finite distributive subsemilattice $\Lambda'$
containing $E$, and $X$ is contained in the submonoid $M'=
\bigsqcup_{e\in \Lambda'} \bigl(\set{e}\times G_e\bigr)$ of $M$.
The temptation is to replace $M$ by $M'$. However, there is no
guarantee that $M'$ satisfies the second part of the MVP, and so
we do not know whether $M'$ is a refinement monoid.
\end{remark}

\begin{remark}\label{Rk:UppBd}
The proof above yields an explicit upper bound for the cardinality
of~$N$ (the desired finite submonoid containing $X$), as a
function of $m$ (fixed positive integer such that $X\subseteq
M[m]$) and $n=|X|$. Now $\bD$ is the sublattice of $\Id\Lambda$
generated by $X\cup\set{0}$. For fixed $x\in X$, we pick elements
$g_{P,x}\in H_{P}$, for $P\in\J_{\bD}([0,e_x])$, such that
$g_x=\sum_{P\in\J_{\bD}([0,e_x])} g_{P,x}$; then put
$U_{P}=\setm{g_{P,x}}{x\in X,\ [0,e_x]\supseteq P}$ and we define
$H'_{P}$ as the subgroup of $H_{P}$ generated by $U_{P}$, for all
$P\in\J(\bD)$. By definition, the subgroups $H'_{P}$ satisfy
\eqref{Eq:H'fplarge}. Hence, the subset
  \[
  Y=\bigcup_{P\in\J(\bD)}\bigl(\set{\varphi(P)}\times U_{P}\bigr)
  \]
is a generating subset of the submonoid $N$ of the proof of
Theorem~\ref{T:FinRep}, with $|Y|\leq|\J(\bD)|\cdot n$. Since
$\bD$ is distributive, every element of $\bD$ is a supremum of
infima of elements of the form $[0,e_x]$, thus every \jirr\
element of $\bD$ has the form $\bigwedge_{x\in I} [0,e_x]$, for
some subset $I$ of $X$. Therefore, $|\J(\bD)|\leq2^n$, and hence,
since $N\subseteq M[m]$, we obtain the estimates
$|N|\leq(m+1)^{|Y|}\leq(m+1)^{2^nn}$.
\end{remark}

We are now ready to establish the key result of the paper, namely
that $\Rep=\cL$.

\begin{theorem}\label{T:Main}
An abelian monoid $M$ is a direct limit of finite direct sums of
monoids of the form $\znzo$ if and only if
\begin{itemize}
\item[(a)] $M$ is a strongly periodic conical refinement monoid.
\item[(b)] For all idempotents $e\leq f$ in $M$, the homomorphism
$\AC{M}{e} \rightarrow \AC{M}{f}$ given by $x\mapsto x+f$ is
injective, and $\AC{M}{e}+f$ is a pure subgroup of $\AC{M}{f}$.
\end{itemize}
\end{theorem}

\begin{proof}
Proposition~\ref{P:NecessLL}, Theorem~\ref{T:FinRep},
Proposition~\ref{P:FinRep}, and Lemma~\ref{L:DirProd}.
\end{proof}

Of course, in case $M$ is countable, the direct limit of
Theorem~\ref{T:Main} may be taken indexed by the natural numbers.

It is easy to restrict the set of cyclic groups used as building
blocks in the theorem, as follows.

\begin{corollary}\label{C:LimRestrPer}
Let $\fm$ be a generalized integer and $M$ an abelian monoid.
Then~$M$ is a direct limit of finite direct sums of monoids of the
form $\znzo$ with $n\mid\fm$ if and only if $M$ satisfies the
conditions of Theorem~\textup{\ref{T:Main}} and
\begin{itemize}
\item[(c)] The order of each element of $M$ divides $\fm$.
\end{itemize}
\end{corollary}

\begin{proof}
We verify the nontrivial direction, $(\Longleftarrow)$. By
Theorem~\ref{T:Main}, $M$ is the direct limit of a direct system
of monoids $M_i$ and transition maps $f_{ij}\colon M_i\rightarrow
M_j$ where each $M_i$ is a finite direct sum of monoids of the
form $\znznz$. It is routine to verify that each $f_{ij}$ maps
$M_i[\fm]$ to $M_j[\fm]$, and that $M[\fm]$ is the direct limit of
the restricted system $\bigl( M_i[\fm],\, f_{ij}|_{M_i[\fm]}
\bigr)$. Assumption (c) says that $M=M[\fm]$, and it only remains
to observe that each $M_i[\fm]$ is a finite direct sum of monoids
$\znznz$ with $n\mid\fm$.
\end{proof}

For the applications to C*-algebras, we need to incorporate
order-units into our direct limits. Recall that an
\emph{order-unit\/} in an abelian monoid $M$ is an element $u\in
M$ such that each $x\in M$ satisfies $x\leq nu$ for some
$n\in\NN$. (In case $M$ is regular, the condition for $u$ to be an
order-unit becomes ``$x\leq u$ for all $x\in M$'', because $2u\leq
u$.) We now work in the category whose objects are pairs $(M,u)$
consisting of abelian monoids $M$ paired with specified
order-units $u$, and whose morphisms are normalized monoid
homomorphisms, that is, a morphism from $(M,u)$ to $(M',u')$ is
any monoid homomorphism from $M$ to $M'$ that sends $u$ to $u'$.
The existence and form of isomorphisms, direct limits, and direct
products in this category are clear. We use the term ``direct
product'' rather than ``direct sum'' here because the natural
construction (via Cartesian products) produces categorical
products which are not coproducts.

Given $m\in\ZZ$ and $n\in\NN$, let us write $\mbar$ for the coset
$m+n\ZZ$, viewed as an element of the monoid $\znznz$; we observe
that $\mbar$ is an order-unit for this monoid.

\begin{corollary}\label{C:Limorderunit}
Let $(M,u)$ be an abelian monoid with order-unit. Then $(M,u)$ is
a direct limit of finite direct products of pairs of the form
$(\znzo,\, \mbar)$ if and only if $M$ satisfies the conditions of
Theorem~\textup{\ref{T:Main}}.
\end{corollary}

\begin{proof} The implication $(\Longrightarrow)$ is immediate from
Theorem~\ref{T:Main}. Conversely, if $M$ satisfies the conditions
of the theorem, then $M$ is the direct limit of a direct system of
monoids $M_i$ and transition maps $f_{ij}$ where each $M_i$ is a
finite direct product of monoids of the form $\znznz$. Let $I$
denote the directed set indexing this direct system, and
$g_i\colon M_i\rightarrow M$ the limiting maps. There exist
$i_0\in I$ and $u_{i_0}\in M_{i_0}$ such that $g_{i_0}(u_{i_0})
=u$. After replacing $I$ by the cofinal subset $\{i\in I\mid i\geq
i_0\}$, we may assume that $i_0$ is the least element of $I$. Set
$u_i= f_{i_0i}(u_{i_0})\in M_i$ for all $i$, so that $g_i(u_i)=
u$.

Next, set $M'_i= \{x\in M_i\mid x\leq u_i\}$ for all $i$, and
observe that $M'_i$ is a submonoid of $M_i$ (remember that
$2u_i\leq u_i$). Moreover, $u_i$ is an order-unit for $M'_i$. Now
any $y\in M$ satisfies $y\leq u$, whence $y= g_i(x)$ for some
$i\in I$ and $x\in M_i$ satisfying $x\leq u_i$, that is, $x\in
M'_i$. Thus, $(M,u)$ is a direct limit of the pairs $(M'_i,u_i)$.
It is straightforward to verify that each $(M'_i,u_i)$ is a finite
direct product of pairs of the form $(\znznz,\, \mbar)$.
\end{proof}

\section{Cuntz limits}
\label{S:C*applics}

Recall that we are using the term \emph{Cuntz limit\/} as an
abbreviation for ``C* inductive limit of a sequence of finite
direct products of full matrix algebras over Cuntz algebras
$\cO_n$ for $n\in\NN$''. (In particular, we are not incorporating
the algebra $\cO_\infty$ into our scheme.) We summarize various
standard facts about the monoids $V(A)$ that will be needed in
applying our monoid-theoretic results to C*-algebras.

First, $V(-)$ is a functor from C*-algebras to abelian monoids
that preserves finite direct products and inductive (direct)
limits \cite[(5.2.3)--(5.2.4)]{Black}. Further, $V(M_m(A))\cong
V(A)$ for any $m\in\NN$ and any $A$, and $V(A)$ is countable if
$A$ is separable \cite[p.~28]{Black}. It is routine to check that
for any unital C*-algebra $A$, the class~$[1_A]$ is an order-unit
in $V(A)$, and that the canonical isomorphism
$V(M_m(A))\rightarrow V(A)$ sends $[1_{M_m(A)}]$ to $m[1_A]$.

The basic K-theoretic information concerning the Cuntz algebras
$\cO_n$ is usually summarized in the statements $K_0(\cO_n)\cong
\ZZ/(n-1)\ZZ$ and $K_1(\cO_n)=0$ \cite[Theorems 3.7--3.8]{Cuntz2}.
However, Cuntz also showed that the Murray-von Neumann equivalence
classes of nonzero projections in $\cO_n$ form a subgroup of
$V(\cO_n)$ which maps isomorphically onto $K_0(\cO_n)$ under the
natural map $V(\cO_n) \rightarrow K_0(\cO_n)$
\cite[p.~188]{Cuntz2}. In addition, the relation $n\cdot 1_{\cO_n}
\sim 1_{\cO_n}$ (a direct consequence of the defining relations
for $\cO_n$) implies that every projection in a matrix algebra
over $\cO_n$ is equivalent to a projection in $\cO_n$ itself. It
follows that $V(\cO_n)\setminus\set{0}$ is a group isomorphic to
$K_0(\cO_n)$, that is, $V(\cO_n)\cong
(\ZZ/(n-1)\ZZ)\sqcup\set{0}$. It is routine to check that this
isomorphism sends $[1_{\cO_n}]$ to the coset $\overline{1}$ in
$\ZZ/(n-1)\ZZ$, and thus we have
  \begin{equation}\label{Eq:VMmOn}
  \bigl( V(M_m(\cO_n)),\, [1_{M_m(\cO_n)}] \bigr) \cong \bigl(
(\ZZ/(n-1)\ZZ)\sqcup\set{0},\, \mbar \bigr)
  \end{equation}
for all $m\geq 1$ and $n\geq 2$. The remaining basic fact that we
shall need is the following lemma. It is essentially equivalent to
\cite[Lemma~6.1]{Ro}; we sketch a proof for the reader's
convenience.

\begin{lemma}\label{L:induceVmaps}
Let $A$ be a finite direct product of full matrix algebras over
Cuntz algebras, $B$ a C*-algebra, and $q\in B$ a projection. Then
any normalized monoid homomorphism
  \[
  \alpha\colon (V(A),[1_A]) \rightarrow (V(B),[q])
  \]
is induced by a C*-algebra map $\phi\colon A\rightarrow B$ that
sends $1_A$ to $q$. That is, $V(\phi)=\alpha$.
\end{lemma}

\begin{proof}
Write $A= \bigoplus_{j=1}^r M_{k_j}(\cO_{n_j})$ for some
$k_j,n_j\in\NN$, and let $p_1,\dots,p_r$ be the corresponding
orthogonal central projections in $A$ summing to $1_A$. Each $p_j$
is an orthogonal sum of pairwise equivalent projections
$e_1^{(j)}$, \dots, $e_{k_j}^{(j)}$ such that $e_1^{(j)}Ae_1^{(j)}
\cong \cO_{n_j}$. In $V(A)$, we have $n_j[e_1^{(j)}]= [e_1^{(j)}]$
for all $j$ and
 \[
 \sum_{j=1}^r k_j[e_1^{(j)}]= \sum_{j=1}^r\, [p_j]= [1_A],
 \]
whence $n_j\alpha([e_1^{(j)}])= \alpha([e_1^{(j)}])$ and
$\sum_{j=1}^r k_j\alpha([e_1^{(j)}])= [q]$ in $V(B)$.
Consequently, $q$ is an orthogonal sum of projections
$q_1,\dots,q_r$ such that $k_j\alpha([e_1^{(j)}])= [q_j]$, and
each~$q_j$ is an orthogonal sum of pairwise equivalent projections
$f_1^{(j)}$, \dots, $f_{k_j}^{(j)}$ such that
$\alpha([e_1^{(j)}])= [f_1^{(j)}]$.

Since $n_j[f_1^{(j)}]= [f_1^{(j)}]$, the projection $f_1^{(j)}$ is
an orthogonal sum of $n_j$ projections each equivalent to
$f_1^{(j)}$, and so there exist $t_1^{(j)}$, \dots, $t_{n_j}^{(j)}
\in f_1^{(j)}Bf_1^{(j)}$ such that $(t_l^{(j)})^*t_m^{(j)}=
\delta_{lm}f_1^{(j)}$ and $\sum_{l=1}^{n_j} t_l^{(j)}(t_l^{(j)})^*
=f_1^{(j)}$. Consequently, there exists a unital C*-algebra map
$\phi_j\colon \cO_{n_j} \rightarrow f_1^{(j)}Bf_1^{(j)}$. Define a
C*-algebra map
 \[
 \phi= \bigoplus_{j=1}^r M_{k_j}(\phi_j)\colon A \longrightarrow
\bigoplus_{j=1}^r M_{k_j}(f_1^{(j)}Bf_1^{(j)}) \cong
\bigoplus_{j=1}^r q_jBq_j \subseteq B.
 \]
It follows from the definition of $\phi$ that $\phi(1_A)=q$ and
$[\phi(e_1^{(j)})]= [f_1^{(j)}]$ for all $j$. Since the classes
$[e_1^{(1)}],\dots,[e_1^{(r)}]$ generate $V(A)$, we conclude that
$V(\phi)=\alpha$.
\end{proof}

\begin{theorem}\label{T:Cstarnonun} An abelian monoid $M$
is isomorphic to $V(A)$ for some Cuntz limit~$A$ if and only if
\begin{itemize}
\item[(a)] $M$ is a countable, strongly periodic, conical refinement
monoid.
\item[(b)] For all idempotents $e\leq f$ in $M$, the homomorphism
$\AC{M}{e} \rightarrow \AC{M}{f}$ given by $x\mapsto x+f$ is
injective, and $\AC{M}{e}+f$ is a pure subgroup of $\AC{M}{f}$.
\end{itemize}
\end{theorem}

\begin{proof} $(\Longrightarrow)$: Recall \eqref{Eq:VMmOn}. Since
$V(-)$ preserves direct limits and finite direct products, the
present implication follows from Theorem~\ref{T:Main}.

$(\Longleftarrow)$: Since $M$ is countable, Theorem~\ref{T:Main}
implies that $M$ is the direct limit of a sequence of the form
\[
  {
  \xymatrix{M_1 \ar[r]^-{\alpha_1} & M_2 \ar[r]^-{\alpha_2}
   & M_3 \ar[r]^-{\alpha_3} & \ar@{}[r]|{\displaystyle\cdots} &
  }}
  \]
where each $M_i$ is a finite direct product of monoids
$\nz{(\ZZ/n_{ij}\ZZ)}$ for some $n_{ij}\in\NN$. Hence, if $A_i$ is
the direct product of the Cuntz algebras $\cO_{n_{ij}+1}$ for the
corresponding indices $j$, then there exists an isomorphism
$h_i\colon V(A_i) \rightarrow M_i$. Each of the homomorphisms
  \[
  h_{i+1}^{-1}\alpha_ih_i\colon V(A_i) \longrightarrow
  V(A_{i+1})
  \]
sends $[1_{A_i}]$ to the class of a projection in $A_{i+1}$, and
so, by Lemma~\ref{L:induceVmaps}, $h_{i+1}^{-1}\alpha_ih_i$ is
induced by a C*-algebra map $\phi_i\colon A_i\rightarrow A_{i+1}$.
Therefore $M \cong V(A)$ where $A$ is the C* inductive limit of
the sequence
\begin{equation*}
  {
  \xymatrix{A_1\ar[r]^-{\phi_1} & A_2\ar[r]^-{\phi_2}
   & A_3\ar[r]^-{\phi_3} & \ar@{}[r]|{\displaystyle\cdots} &
  }}\tag*{\qed}
  \end{equation*}
\renewcommand{\qed}{}
\end{proof}

A structural description of the  monoids appearing in
Theorem~\ref{T:Cstarnonun} is easily obtained with the help of
Theorem~\ref{T:RAMstruct}, as follows.

\begin{corollary}\label{C:VAnonun}
Let $M$ be an abelian monoid. Then $M\cong V(A)$ for some Cuntz
limit $A$ if and only if
  \[
  M \cong \bigsqcup_{e\in\Lambda} \bigl( \set{e}\times G_e
\bigr) \subseteq \Lambda\times G
  \]
where
\begin{itemize}
\item[(a)] $\Lambda$ is a countable distributive semilattice.
\item[(b)] $G$ is a countable torsion abelian group.
\item[(c)] $G_e$ is a pure subgroup of $G$ for all $e\in\Lambda$.
\item[(d)] $G_0=\set{0}$ and $\bigcup_{e\in\Lambda} G_e =G$.
\item[(e)] $G_e+G_f= G_{e+f}$ and $G_e\cap G_f=
\bigcup_{g\in\Lambda,\, g\leq e,f} G_g$ for all $e$,
$f\in\Lambda$.
\end{itemize}
\end{corollary}

We can also characterize the monoids $V(A)$ for Cuntz limits~$A$
with a restricted set of building blocks $\cO_n$, as follows.

\begin{corollary}\label{C:Cstarnonun}
Let $M$ be an abelian monoid and $\fm$ a generalized integer. Then
$M\cong V(A)$ for some C* inductive limit of a sequence of finite
direct products of full matrix algebras over Cuntz algebras
$\cO_n$ with $n-1\mid\fm$ if and only if $M$ satisfies the
conditions of Theorem~\textup{\ref{T:Cstarnonun}} and the order of
each element of $M$ divides $\fm$.
\end{corollary}

\begin{proof}
Theorem~\ref{T:Cstarnonun} and Corollary~\ref{C:LimRestrPer}.
\end{proof}

Finally, we establish the unital cases of the above results.

\begin{theorem}\label{T:seastar} Let $(M,u)$ be an abelian monoid with
order-unit. Then $(M,u)\cong (V(A),[1_A])$ for some unital Cuntz
limit $A$ if and only if
\begin{itemize}
\item[(a)] $M$ is a countable, strongly periodic, conical refinement
monoid.
\item[(b)] For all idempotents $e\leq f$ in $M$, the homomorphism
$\AC{M}{e} \rightarrow \AC{M}{f}$ given by $x\mapsto x+f$ is
injective, and $\AC{M}{e}+f$ is a pure subgroup of $\AC{M}{f}$.
\end{itemize}
\end{theorem}

\begin{proof} $(\Longrightarrow)$: Theorem~\ref{T:Cstarnonun}.

$(\Longleftarrow)$: Corollary~\ref{C:Limorderunit} implies that
$(M,u)$ is the direct limit of a sequence of the form
\[
  {
  \xymatrix{(M_1,u_1)\ar[r]^-{\alpha_1} &
(M_2,u_2)\ar[r]^-{\alpha_2}
   & (M_3,u_3)\ar[r]^-{\alpha_3} & \ar@{}[r]|{\displaystyle\cdots} &
  }}
  \]
where each $(M_i,u_i)$ is a finite direct product of pairs
$(\nz{(\ZZ/n_{ij}\ZZ)},\, \mbar_{ij})$ for some
$n_{ij},m_{ij}\in\NN$. In view of~\eqref{Eq:VMmOn}, there exist
isomorphisms $h_i\colon (V(A_i),[1_{A_i}]) \rightarrow (M_i,u_i)$
where $A_i$ is the direct product of the matrix algebras
$M_{m_{ij}}(\mathcal{O}_{n_{ij}+1})$. Each of the normalized
homomorphisms
  \[
  h_{i+1}^{-1}\alpha_ih_i\colon (V(A_i),[1_{A_i}]) \longrightarrow
  (V(A_{i+1}),[1_{A_{i+1}}])
  \]
is induced by a unital C*-algebra map $\phi_i\colon A_i\rightarrow
A_{i+1}$ (Lemma~\ref{L:induceVmaps}). Therefore $(M,u) \cong
(V(A),[1_A])$ where $A$ is the C* inductive limit of the sequence
\begin{equation*}
  {
  \xymatrix{A_1\ar[r]^-{\phi_1} & A_2\ar[r]^-{\phi_2}
   & A_3\ar[r]^-{\phi_3} & \ar@{}[r]|{\displaystyle\cdots} &
  }}\tag*{\qed}
  \end{equation*}
\renewcommand{\qed}{}
\end{proof}

\begin{corollary}\label{C:VAstruct}
Let $(M,u)$ be an abelian monoid with order-unit. Then $(M,u)\cong
(V(A),[1_A])$ for some unital Cuntz limit $A$ if and only if
  \[
  (M,u) \cong \Bigl( \bigsqcup_{e\in\Lambda} \bigl( \set{e}\times G_e
\bigr),\, (1,u_1) \Bigr) \subseteq \bigl( \Lambda\times G_1,\,
(1,u_1) \bigr)
 \]
where
\begin{itemize}
\item[(a)] $\Lambda$ is a countable distributive semilattice with
maximum element $1$.
\item[(b)] $G_1$ is a countable torsion abelian group.
\item[(c)] $G_e$ is a pure subgroup of $G_1$ for all $e\in\Lambda$,
and $G_0=\set{0}$.
\item[(d)] $G_e+G_f= G_{e+f}$ and $G_e\cap G_f=
\bigcup_{g\in\Lambda,\, g\leq e,f} G_g$ for all $e$,
$f\in\Lambda$.
\item[(e)] $u_1\in G_1$.
\end{itemize}
\end{corollary}

\begin{proof}
$(\Longrightarrow$): By Corollary~\ref{C:VAnonun}, $M$ is
isomorphic to a monoid of the form
  \[
  M'= \bigsqcup_{e\in\Lambda} \bigl( \set{e}\times G_e
 \bigr) \subseteq \Lambda\times G
  \]
for some countable distributive semilattice $\Lambda$ and some
countable torsion abelian group $G$ with subgroups $G_e$
satisfying the conditions of that corollary. An isomorphism
$M\rightarrow M'$ must carry $u$ to an order-unit $u'=
(\eps,u_\eps) \in M'$. For each $e\in\Lambda$, there exists
$n\in\NN$ such that $(e,0)\leq nu'= (\eps,nu_\eps)$, whence
$e\le\eps$. Thus, $\eps$ is the largest element of $\Lambda$, and
we rename it in the standard way: $\eps=1$. Conditions (a)--(e)
are now all satisfied.

$(\Longleftarrow)$: With the help of Theorem~\ref{T:RAMstruct}, it
is clear that $M$ satisfies conditions (a) and (b) of
Theorem~\ref{T:seastar}.
\end{proof}

\begin{corollary}\label{C:seastar}
Let $(M,u)$ be an abelian monoid with order-unit, and $\fm$ a
generalized integer. Then $(M,u)\cong (V(A),[1_A])$ for some
unital C* inductive limit of a sequence of finite direct products
of full matrix algebras over Cuntz algebras $\cO_n$ with
$n-1\mid\fm$ if and only if $M$ satisfies the conditions of
Theorem~\textup{\ref{T:seastar}} and the order of each element of
$M$ divides $\fm$.
\end{corollary}

\begin{proof}
Theorem~\ref{T:seastar} and Corollary~\ref{C:LimRestrPer}.
\end{proof}

\section*{Acknowledgments}

Part of this work was done during visits of the second author to
the Department of Mathematics of the University of California at
Santa Barbara (USA) and the D\'epartement de Math\'ematiques de
l'Universit\'e de Caen (France). The second author wants to thank
both host centers for their warm hospitality.

\end{document}